\documentclass[12pt]{amsart}
\usepackage{amsfonts}
\usepackage{mathrsfs}
\usepackage{graphicx}
\usepackage{amsmath}
\usepackage{amssymb}
\usepackage{amsthm}
\usepackage{slashbox}
\usepackage{multirow}
\usepackage[top=1in,bottom=1in,left=1in,right=1in]{geometry}

\newtheorem{theorem}{Theorem}

\def\Ai{\mathrm {Ai}}

\renewcommand\O[1]{\mathcal{O}\left({#1}\right)}

\begin{document}

\title[Global asymptotics of Stieltjes-Wigert polynomials]
{Global asymptotics of Stieltjes-Wigert polynomials}

\author{Y. T. Li}
\email{yutianli@hkbu.edu.hk}
\address{Institute of Computational and Theoretical Studies, and Department of Mathematics,
Hong Kong Baptist University, Kowloon, Hong Kong}
\thanks{The work of Y. T. Li is supported by the HKBU Strategic Development Fund and a fund from HKBU (no.: 38-40-106)}

\author{R. Wong}
\email{mawong@cityu.edu.hk}
\address{Liu Bie Ju Centre for Mathematical Sciences, City University of Hong Kong, Kowloon, Hong Kong}

\maketitle

\begin{abstract}
Asymptotic formulas are derived for the
Stieltjes-Wigert polynomials $S_n(z;q)$ in the complex plane as the degree $n$ grows to infinity.
One formula holds in any disc centered at the origin, and the other holds outside any smaller disc centered at the origin;
the two regions together cover the whole plane.
In each region, the $q$-Airy function $A_q(z)$ is used as the approximant.
For real $x>1/4$,
a limiting relation is also established between the $q$-Airy function $A_q(x)$ and the ordinary Airy function $\Ai(x)$ as $q\to1$.
\end{abstract}

\section{introduction}
We first fix some notations. Let $k>0$ be a fixed number and
\begin{equation}\label{eq:1}
q=\exp\{-(2k^2)^{-1}\}.
\end{equation}
Note that $0<q<1$. The $q$-shifted factorial is given by
$$
(a,q)_0=1,\qquad
(a,q)_n=\prod_{j=0}^{n-1}(1-aq^j),\qquad n=1,2,\cdots.
$$
With these notations, the Stieltjes-Wigert polynomials
\begin{equation}\label{eq:1.4}
S_n(x;q)=\sum_{j=0}^n\frac{q^{j^2}}{(q,q)_j(q;q)_{n-j}}(-x)^j,\qquad n=0,1,2,\cdots,
\end{equation}
are orthogonal with respect to the weight function
\begin{equation}\label{eq:1.5}
w(x)=k\pi^{-\frac12}\exp\{-k^2\log^2x\}
\end{equation}
for $0<x<\infty$; see \cite[(18.27.18)]{NIST} and \cite[(3.27.1)]{KS}.
It should be mentioned that the Stieltjes-Wigert polynomials belong to the
indetermined moment class and the weight function in (\ref{eq:1.5}) is not unique; see~\cite{Christiansen}.
One important property of the Stieltjes-Wigert polynomials is the symmetry relation
\begin{equation}\label{eq:5}
S_n(z;q)=(-zq^n)^n S_n\left(\frac{1}{zq^{2n}};q\right),
\end{equation}
which can be easily verified by changing the index $j$ to $n-j$ in the explicit expression given in (\ref{eq:1.4}).
In some literatures, the variable $x$ in (\ref{eq:1.4}) is replaced by $q^{\frac12}x$; see, for example, Szeg\"o~\cite{Szego}, Chihara~\cite{Chihara},
and Wang and Wong~\cite{WangWong}. The notation for the Stietjes-Wigert polynomials used in these literatures is
\begin{equation}\label{eq:SWp}
p_n(x)=(-1)^nq^{n/2+1/4}\sqrt{(q;q)_n}S_n(q^\frac12 x;q).
\end{equation}

The Stietjes-Wigert polynomials appear in random walks and random matrix formulation of Chern-Simons theory on Seifert
manifolds; see \cite{BaikSuidan,DT}.

The asymptotics of the Stieltjes-Wigert polynomials, as the degree tends to infinity, has been studied by several authors.
In 1923, Wigert~\cite{Wigert} proved that the polynomials have the limiting behavior
\begin{equation}
\lim_{n\to\infty}(-1)^nq^{-n/2}p_n(x)=\frac{q^{1/4}}{\sqrt{(q;q)_\infty}}\sum_{k=0}^\infty
(-1)^k\frac{q^{k^2+k/2}}{(q;q)_k}x^k,
\end{equation}
which can be put in terms of the $q$-Airy function
\begin{equation}\label{eq:8}
A_q(z)=\sum_{k=0}^\infty\frac{q^{k^2}}{(q,q)_k}(-z)^k.
\end{equation}
This function appeared in the third identity on p.57 of Ramanujan's ``Lost Notebook''~\cite{Ramanujan}.
(For this reason, it is also known as the Ramanujan function.)
In terms of the $q$-Airy function, Wigert's result can be stated as
\begin{equation}\label{eq:9}
\lim_{n\to\infty}S_n(x;q)=\frac{1}{(q;q)_\infty}A_q(x).
\end{equation}

It is known that all zeros of $S_n(x;q)$ lie in the interval $(0,4q^{-2n})$; see~\cite{WangWong}.
Hence, we introduce a new scale
\begin{equation}
z:=q^{-nt}u
\end{equation}
with $u\in\mathbb C\setminus\{0\}$ and $t\in\mathbb R$.
The values of $t=0$ and $t=2$ can be regarded as the turning points of $S_n(q^{-nt}u;q)$.
Taking into account the symmetry relation in (\ref{eq:5}),
one may restrict oneself to the case $t\geq 1$; see~\cite[(1.4)]{WangXSWong2}.
(However, in the present paper, we will not make this restriction.)
The case $t=2$ has been studied by Ismail~\cite{Ismail2},
and he proved
\begin{equation}
\lim_{n\to\infty}q^{n^2(t-1)}(-u)^{-n}S_n(uq^{-nt};q)=\frac{1}{(q;q)_\infty}A_q\left(\frac{q^{n(t-2)}}{u}\right),\qquad t=2,
\end{equation}
uniformly on compact subsects of $\mathbb C\setminus\{0\}$; see \cite[Theorem 2.5]{Ismail2}.
This result can in fact be derived directly from Wigert's result in (\ref{eq:9}) via the symmetry relation mentioned in (\ref{eq:5}).
In~\cite{IsmailZhang}, Ismail and Zhang extended the validity of this result to $t\geq2$.
For $1\leq t<2$, Ismail and Zhang~\cite{IsmailZhang} gave asymptotic formulas for these polynomials in terms of
the theta-type function
\begin{equation}
\Theta_q(z)=\sum_{k=-\infty}^\infty q^{k^2}z^k,
\end{equation}
but in a very complicated manner.
The result in~\cite{IsmailZhang} was then simplified by Wang and Wong~\cite{WangXSWong1}.
For instance, when $1\leq t<2$, Wang and Wong proved that
\begin{equation}
S_n(uq^{-nt};q)=\frac{(-u)^{n-m}q^{n^2(1-t)-m[n(2-t)-m]}}{(q;q)_n(q;q)_\infty}
\left\{\Theta_q\left(\frac{q^{2m-n(2-t)}}{-u}\right)+\O{q^{n(l-\delta)}}\right\},
\end{equation}
where $l=\frac12(2-t)$, $m=\lfloor nl\rfloor$ and $\delta>0$ is any small number; see~\cite[Corollary 2]{WangXSWong1}.
Note that all these results are not valid in a neighborhood containing $t=2$, one of the turning points.
To resolve this issue, a uniform asymptotic formula was given by Wang and Wong in a second paper~\cite{WangXSWong2}.
For $z:=uq^{-nt}$ with $t>2(1-\delta)$, $\delta$ being any small positive constant, they showed that
\begin{equation}\label{eq:14}
S_n(z;q)=\frac{(-z)^nq^{n^2}}{(q;q)_n}\big[A_{q,n}(q^{-2n}/z)+r_n(z)\big],
\end{equation}
where $r_n(z)$ is the remainder and $A_{q,n}(z)$ is the $q$-Airy polynomial obtained by truncating the infinite series in (\ref{eq:8}) at $k=n$, \textit{i.e.},
$$
A_{q,n}(z)=\sum_{k=0}^n \frac{q^{k^2}}{(q,q)_k}(-z)^k.
$$

In this paper, we shall show that the $q$-Airy polynomial $A_{q,n}(z)$ in equation (\ref{eq:14}) can be replaced by the $q$-Airy function $A_q(z)$.
Moreover, we shall show that the resulting formula is global.
More precisely, we have the following result.

\begin{theorem}\label{thm:1}
Let $z:=uq^{-nt}$ with $-\infty<t<2$, $u\in\mathbb C$ and $|u|\leq R$, where $R>0$ is any fixed positive number. We have
\begin{equation}\label{eq:3.13}
S_{n}(z;q)=\frac{1}{(q;q)_n}\left[A_q(z)+r_n(z)\right],
\end{equation}
where the remainder satisfies
\begin{equation}\label{eq:3.14}
|r_n(z)|\leq \left[\frac{q^{n(1-\sigma)}}{1-q}+\frac{2}{1-q}\left(\frac{1}{2}\right)^{\lfloor n\sigma\rfloor}\right]A_q(-|z|)
\end{equation}
with $\sigma=\max\{\frac12,\frac12+\frac t4\}$.

Let $z:=uq^{-nt}$ with $0<t<\infty$, $u\in\mathbb C$ and $|u|\geq 1/R$, where $R>0$ is any fixed positive number. We have
\begin{equation}\label{eq:3.15}
S_{n}(z;q)=\frac{(-z)^nq^{n^2}}{(q;q)_n}\left[A_q(q^{-2n}/z)+r_n(z)\right],
\end{equation}
where the remainder satisfies
\begin{equation}\label{eq:3.16}
|r_n(z)|\leq \left[\frac{q^{n(\delta-1)}}{1-q}+\frac{2}{1-q}\left(\frac{1}{2}\right)^{\lfloor n(2-\delta)\rfloor}\right]A_q\big(-q^{-2n}/|z|\big)
\end{equation}
with $\delta=\min\{\frac32,1+\frac t4\}$.
\end{theorem}

Note that the quantities inside the square brackets in (\ref{eq:3.14}) and (\ref{eq:3.16}) are exponentially small.
Furthermore, the sizes of the $q$-Airy functions in these two equations
for large values of their arguments are about the same as
the leading terms in their corresponding approximation formulas
in (\ref{eq:3.13}) and (\ref{eq:3.15}); cf. (\ref{eq:lim2}) below.

The investigations mentioned above all started with the explicit expression of $S_n(x;q)$ given in (\ref{eq:1.4}).
In~\cite{WangWong}, Wang and Wong used a different method, namely, the Riemann-Hilbert approach,
to get a uniform asymptotic expansion of the Stietjes-Wigert polynomials $p_n(x)$
in terms of Airy functions.
But, the main result in \cite{WangWong} needs a correction,
and the correction is that the parameter $k$ in (\ref{eq:1}) and (\ref{eq:1.5})
should depend on $n$ and tend to infinity as $n$ tends to infinity.
In other words, the result in~\cite{WangWong} holds with a varying weight and
the number $q$ in (\ref{eq:1}) is required to approach 1.
In fact, more precisely, if $k\sim n^\sigma$ as $n\to\infty$ and $0<\sigma<\frac12$,
then all formulas in \cite{WangWong} remain valid;
if $\sigma\geq\frac12$, then some equations need be amended,
but the main result still holds;
see also Baik and Suidan~\cite{BaikSuidan}.
Note that the results in Theorem~\ref{thm:1} hold even when $q$ tends to 1.
Comparing the results in Theorem~\ref{thm:1} and in \cite{WangWong},
we will establish the following limiting relation between the $q$-Airy function and the ordinary Airy function:
\begin{theorem}\label{thm:2}
Let $\xi(x)$ denote the function defined by
\begin{equation}\label{eq:12}
\frac23[\xi(x)]^{3/2}=\frac{1}{\log(1/q)}\int_0^{\log(4x)}\arctan\sqrt{e^s-1}ds.
\end{equation}
Then, for fixed $x>1/4$, the $q$-Airy function has the following asymptotic approximation
\begin{equation}\label{eq:13}
A_q(\sqrt q x)\sim 2\sqrt{\pi }\exp\left\{\frac{3\log^2x-\pi^2}{12\log(1/q)} \right\}\left(\frac{\xi(x)}{4x-1}\right)^\frac14 \Ai(-\xi(x))
\end{equation}
as $q\to1$.
\end{theorem}
This result in fact holds for any $x\in\mathbb C\setminus\{0\}$, if we replace $\xi(x)$ in (\ref{eq:13}) by $\widetilde{\xi}(x)$ defined by
\begin{equation*}
\frac23\big[-\widetilde{\xi}(x)\big]^{3/2}=\frac{1}{\log q}\left[\int_0^{\log(4x)}\log\left(1+\sqrt{1-e^s}\right)ds-\frac{\log^2(4x)}{4}\right].
\end{equation*}
In view of this result, $A_q$ is indeed a $q$-analogue of the Airy function.

\section{Global asymptotics of Stieltjes-Wigert polynomials}

\textbf{Proof of Theorem 1.}

If $u=0$, then $z=uq^{-nt}=0$ and
$$
(q;q)_nS_n(0;q)=A_q(0)=1
$$
by the series representations in (\ref{eq:1.4}) and (\ref{eq:8}).
Hence, the result in (\ref{eq:3.13})-(\ref{eq:3.14}) follows immediately,
and we just need consider the case when $u\neq0$.
Noting that
$$
z=uq^{-nt}=\frac{u}{|u|}q^{-n\big(t-\frac{\log|u|}{n\log q}\big)},
$$
we may assume $|u|=1$ without loss of generality.

For notational convenience, we put
\begin{equation}\label{eq:3.17}
\begin{aligned}
r_n(z)=&(q;q)_n S_n(z;q)-A_q(z)\\
=&\sum_{j=0}^n\left[\frac{(q;q)_n}{(q;q)_{n-j}}-1\right]\frac{q^{j^2}}{(q;q)_j}(-z)^j
-\sum_{j=n+1}^\infty \frac{q^{j^2}}{(q;q)_j}(-z)^j\\
=:&I_1+I_2+I_3,
\end{aligned}
\end{equation}
where
$$
I_1=\sum_{j=0}^{\lfloor n\sigma\rfloor}
\left[\frac{(q;q)_n}{(q;q)_{n-j}}-1\right]\frac{q^{j^2}}{(q;q)_j}(-z)^j,
$$
$$
I_2=\sum_{j=\lfloor n\sigma\rfloor+1}^{n}
\left[\frac{(q;q)_n}{(q;q)_{n-j}}-1\right]\frac{q^{j^2}}{(q;q)_j}(-z)^j,
$$
$$
I_3=
\sum_{j=n+1}^\infty \frac{q^{j^2}}{(q;q)_j}(-z)^j,
$$
and $0<\sigma<1$ is a constant to be specified later.
In view of the inequality $1-ab<(1-a)+(1-b)$ for any $a,b\in(0,1)$, we have
\begin{equation}
1-\frac{(q;q)_n}{(q;q)_{n-j}}
=1-(q^{n-j+1};q)_j
<\sum_{i=1}^jq^{n-j+i}
<\frac{q^{n-j}}{1-q}.
\label{267488156}
\end{equation}
Thus, we have
\begin{equation}\label{eq:3.18}
|I_1|\leq\sum_{j=0}^{\lfloor n\sigma\rfloor}
\frac{q^{n-j}}{1-q}\frac{q^{j^2}}{(q;q)_j}|z|^j
<\frac{q^{n(1-\sigma)}}{1-q}A_q(-|z|).
\end{equation}

Let
\begin{equation}
\sigma=\max\left\{\frac12,\frac{t}{4}+\frac12\right\},\qquad
l=\max\left\{\frac n4,\lfloor\frac{nt}{2}\rfloor+1\right\}.
\end{equation}
It is readily seen that
$\frac t2<\sigma<1$ and that $\frac {nt}2<l<n$ for sufficiently large $n$.
For $j\geq \lfloor n\sigma\rfloor$ and $n$ sufficiently large, there exists $\gamma>0$ such that
$$
(j^2-ntj)-(l^2-ntl)=(j-l)(j+l-nt)>(j-l)^2\geq \gamma^2 j^2.
$$
Then, it follows that
$$
\begin{aligned}
|I_2+I_3|\leq&
\frac{1}{1-q}\sum_{j=\lfloor n\sigma\rfloor+1}^{\infty}
\frac{q^{j^2}}{(q;q)_j}|z|^j\\
=& \frac{q^{l^2-ntl}}{1-q}\sum_{j=\lfloor n\sigma\rfloor+1}^{\infty}
\frac{q^{(j^2-ntj)-(l^2-ntl)}}{(q;q)_j}\\
<& \frac{q^{l^2-ntl}}{(1-q)(q;q)_\infty}\sum_{j=\lfloor n\sigma\rfloor+1}^{\infty}
q^{\gamma^2 j^2}.
\end{aligned}
$$
Since
$q^{\gamma^2j}\leq\frac12$
for $j\geq \lfloor n\sigma\rfloor$ and sufficiently large $n$, we have
\begin{equation}\label{eq:32145}
\begin{aligned}
|I_2+I_3|
\leq& \frac{q^{l^2-ntl}}{(1-q)(q;q)_\infty}\sum_{j=\lfloor n\sigma\rfloor+1}^{\infty}
\left(\frac12\right)^j\\
\leq& \frac{q^{l^2-ntl}}{(1-q)(q;q)_\infty}
\left(\frac12\right)^{\lfloor n\sigma\rfloor}.
\end{aligned}
\end{equation}
Similar to the inequality in (\ref{267488156}), we have
$$
1-\frac{(q;q)_\infty}{(q;q)_l}=1-(q^{l+1};q)_\infty<\frac{q^l}{1-q}\leq \frac{q^\frac{n}{4}}{1-q}<\frac12
$$
for sufficiently large $n$, which leads to
\begin{equation}\label{eq:2458}
\frac{(q;q)_l}{(q;q)_\infty}<\frac{1}{1-\frac{1}{2}}<2.
\end{equation}
Moreover, we have
\begin{equation}\label{eq:125482}
\frac{q^{l^2-ntl}}{(q;q)_l}\leq \sum_{j=0}^\infty\frac{q^{j^2-ntj}}{(q;q)_j}
\leq \sum_{j=0}^\infty\frac{q^{j^2}}{(q;q)_j}|q^{-nt}u|^{j}
=A_q(-|z|).
\end{equation}
A combination of the inequalities in (\ref{eq:32145}), (\ref{eq:2458}) and (\ref{eq:125482}) gives
$$
|I_2+I_3|\leq \frac{2}{1-q}
\left(\frac12\right)^{\lfloor n\sigma\rfloor}A_q(-|z|),
$$
which, together with (\ref{eq:3.17}) and (\ref{eq:3.18}), yields (\ref{eq:3.14}).

The result in (\ref{eq:3.15})-(\ref{eq:3.16}) can be proved in a similar manner.
One can also obtain this directly from (\ref{eq:3.13}), (\ref{eq:3.14})
and the symmetry relation of $S_n(z;q)$ mentioned in (\ref{eq:5}).

\section{Limiting behavior of $A_q(z)$ as $q\to1$}

It is known that $A_q(z)$ has infinitely many positive zeros and satisfies the three-term
recurrence relation \cite{Ismail}
\begin{equation}
A_q(z)-A_q(qz)+qzA_q(q^2z)=0.
\end{equation}
Moreover, Zhang \cite{Zhang} has shown that
\begin{equation}
\lim_{q\to 1^-}A_q\big((1-q)z\big)=e^{-z}
\end{equation}
for any fixed $z\in\mathbb C$. In \cite[Propsition 1]{WangXSWong2}, Wang and Wong proved that
\begin{equation}\label{eq:lim2}
 A_q(z)=\frac{(-z)^mq^{m^2}}{(q;q)_\infty}\big[\Theta_q(-q^{2m}z)+\O{q^{m(1-\delta)}}\big]
\end{equation}
as $z\to\infty$, where $m:=\lfloor\frac{\ln|z|}{-2\ln q}\rfloor$ and $\delta>0$ is any small number; see also Zhang \cite[Theorem 2.1]{Zhang}.

In this section, we shall establish the limiting relation of $A_q(z)$,
as $q\to 1$, stated in Theorem~\ref{thm:2}.
Let us first review some of the results given in Wang and Wong~\cite{WangWong}.
Let
\begin{equation}\label{eq:32}
k:=k(n)=n^{\frac14}.
\end{equation}
Then
\begin{equation}\label{eq:33}
q=\exp\{-(2k^2)^{-1}\}=\exp\{-(2\sqrt n)^{-1}\}
= 1-\frac{1}{2\sqrt{n}}+\O{\frac1n}
\end{equation}
as $n\to\infty$.
In Sec. 1 we have commented that the main result in Wang and Wong~\cite{WangWong} holds
if the parameter $k$ in (\ref{eq:1}) depends on $n$ and satisfies a growth condition such as the one given in (\ref{eq:32}).
The MRS numbers $\alpha_n$ and $\beta_n$ have been calculated in~\cite{WangWong},
and are explicitly given in equations (3.19) and (3.20).
We note that
\begin{equation}\label{eq:34}
\alpha_n\sim \frac14,\qquad \beta_n\sim 4q^{-(2n+1)},\qquad \text{as } n\to\infty;
\end{equation}
cf. \cite[eq. (2.4)]{WangWong}.
Furthermore, the change of variable $t \to y$ defined by
\begin{equation}\label{eq:t}
y=\sqrt{\alpha_n\beta_n}\exp\left[\frac t2 \log(\beta_n/\alpha_n)\right]
\end{equation}
takes the interval $-1\leq t\leq 1$ onto the interval $\alpha_n\leq y\leq \beta_n$,
where $y=1/(xq^{2n+1})$; see~\cite[(2.14)]{WangWong}.
Let $\pi_n(x)$ denote the monic Stieltjes-Wigert polynomial
$$
\pi_n(x)=\frac{p_n(x)}{\gamma_n},
$$
where $p_n(x)$ is the polynomial given in (\ref{eq:SWp}) and
$$
\gamma_n=q^{n^2+n+1/4}/\sqrt{(q;q)_n};
$$
that is
\begin{equation}
\pi_n(x)=(-1)^nq^{-n^2-n/2}(q;q)_nS_n(q^{1/2}x;q).
\end{equation}
The major results in \cite{WangWong} is the asymptotic formula
\begin{equation}\label{eq:asymp}
\pi_n\big(y\big)=\frac{\sqrt{\pi}e^{l_n/2}}{\sqrt{w(y)}}\left\{N^\frac16 \Ai\left(N^\frac23\eta_n(t)\right)A(y,n)+\O{N^{-\frac16}}\right\},
\end{equation}
where $N=n+\frac12$,
\begin{equation}\label{eq:A}
A(y,n)=\frac{[\eta_n(t)]^{1/4}(\beta_n-\alpha_n)^{1/2}}{[(y-\alpha_n)(y-\beta_n)]^{1/4}},
\end{equation}
$$
\frac23[-\eta_n(t)]^{3/2}=\frac{a}{N\log(1/q)}\int_t^1\arctan \frac{\sqrt{(e^{a\tau}-e^{-a})(e^a-e^{a\tau})}}{e^{a\tau}+1}d\tau
$$
and
$$
a:=\frac12\log\frac{\beta_n}{\alpha_n};
$$
see \cite[(2.11), (2.12) and (6.7)]{WangWong}.
Recall that $y$ in (\ref{eq:asymp}) and (\ref{eq:A}) is a function of $t$, given in (\ref{eq:t}).
Using the formula \cite[(6.8)]{WangWong}
$$
\frac23[-\eta_n(t)]^{3/2}=\frac{1}{N\log(1/q)}\int_0^{a(1-t)}\arctan\sqrt{e^s-1}ds+\O{q^{\frac12\delta N}},
$$
where $-1+\delta<t<1$, we have
\begin{equation}
N^{\frac23}\eta_n\sim-\xi(x)
\label{2359824}
\end{equation}
as $n\to\infty$, where $\xi(x)$ is the function defined in (\ref{eq:12}).
The last equation gives
\begin{equation}
\Ai(N^{2/3}\eta_n)\sim\Ai(-\xi(x)),\qquad n\to\infty.
\label{2985245}
\end{equation}

Next, we show that
\begin{equation}
 N^\frac16 A\left(y,n\right)\sim 2\sqrt{x}\left(\frac{\xi(x)}{4x-1}\right)^\frac14,\qquad n\to\infty.
\label{298888}
\end{equation}
To this end, we note from (\ref{eq:A}) that
\begin{equation}
\begin{aligned}
 N^\frac16 A\left(\frac{1}{xq^{2n+1}},n\right)
=&
\left\{\frac{N^{2/3}\eta_n(t)(\beta_n-\alpha_n)^2}{[1/xq^{2n+1}-\alpha_n][1/xq^{2n+1}-\beta_n]}\right\}^\frac14\\
=&
\left\{\frac{N^{2/3}\eta_n(t)(\beta_n-\alpha_n)^2x^2q^{2(2n+1)}}{(1-\alpha_nxq^{2n+1})(1-\beta_nxq^{2n+1})}\right\}^\frac14.
\end{aligned}
\end{equation}
By using equations (3.17) and (3.18) in~\cite{WangWong}, we have
\begin{equation}\label{eq:14253}
\alpha_n\beta_n=q^{-(2n+1)}.
\end{equation}
Since $\beta_n$ is large,
the last two equations, together with (\ref{2359824}), give
\begin{equation}
N^\frac16 A\left(\frac{1}{xq^{2n+1}},n\right)\sim \left(\frac{-\xi(x)x^2/\alpha_n^2}{1-x/\alpha_n}\right)^\frac14,\qquad n\to\infty,
\end{equation}
thus proving (\ref{298888}).
Here, we have also made use of the fact that $\alpha_n\sim\frac14$ as $n\to\infty$.

Finally, we evaluate the asymptotics of
\begin{equation}
\frac{\sqrt\pi e^{\frac12 l_n}}{y^n\sqrt{w(y)}}=\sqrt{\pi}\exp\left\{\frac12 l_n+\frac12 k^2 \log^2 y-\frac12\log\frac k{\sqrt{\pi}}-n\log y \right\}.
\end{equation}
Recall $y=1/xq^{2n+1}$. Using the formula \cite[(3.29)]{WangWong}
$$
l_n=\frac{N(N-1)}{k^2}-\frac{k^2\pi^2}{3}+\log\frac{k}{\sqrt\pi}+\O{Nq^N},
$$
we obtain
\begin{equation}
\frac{\sqrt\pi e^{\frac12 l_n}}{y^n\sqrt{w(y)}}\sim \sqrt{\pi}\exp\left\{ \frac{N(N-1)}{2k^2}-\frac{k^2\pi^2}{6}+\frac12 k^2 \log^2 (xq^{2n+1})+n\log(xq^{2n+1}) \right\}.
\label{abcdfe}
\end{equation}
Since $\log(1/q)=1/2k^2$ and $N=n+\frac12$, one can show that the right-hand side of (\ref{abcdfe}) is equal to
\begin{equation}
 \sqrt{\pi}\exp\left\{-\frac12\log x+\frac{3\log^2x-\pi^2}{12\log(1/q)} \right\}=\sqrt\frac{\pi}{x}\exp\left\{\frac{3\log^2x-\pi^2}{12\log(1/q)} \right\}.
\end{equation}
Hence, we obtain
\begin{equation}
\frac{\sqrt\pi e^{\frac12 l_n}}{y^n\sqrt{w(y)}}\sim \sqrt\frac{\pi}{x}\exp\left\{\frac{3\log^2x-\pi^2}{12\log(1/q)} \right\},\qquad n\to\infty.
\label{3122334}
\end{equation}
A combination of (\ref{eq:asymp}), (\ref{298888}) and (\ref{3122334}) yields
\begin{equation}
(xq^{2n+1})^n\pi_n(1/xq^{2n+1})\sim \sqrt{\frac\pi x}\exp\left\{\frac{3\log^2x-\pi^2}{12\log(1/q)} \right\}2\sqrt{x}\left(\frac{\xi(x)}{4x-1}\right)^\frac14 \Ai(-\xi(x)).
\label{3288546}
\end{equation}
Therefore, we have
$$
(xq^{2n+1})^n\pi_n(1/xq^{2n+1})=(q;q)_n(-q^{1/2}x q^{n})^n S_n(q^{-2n}/q^{1/2}x;q)=(q;q)_nS_n(q^\frac12 x;q).
$$
Here, we have made use of the symmetry relation of $S_n(z;q)$ given in (\ref{eq:5}).
Note that the result in Theorem~\ref{thm:1} in fact holds with $q$ satisfying (\ref{eq:33}).
Thus, Theorem~\ref{thm:1} gives
\begin{equation}
(xq^{2n+1})^n\pi_n(1/xq^{2n+1})=A_q(\sqrt q x)+r_n(\sqrt q x).
\label{334567}
\end{equation}
Coupling (\ref{3288546}) and (\ref{334567}), we obtain the desired formula
$$
A_q(\sqrt q x)\sim 2\sqrt{\pi }\exp\left\{\frac{3\log^2x-\pi^2}{12\log(1/q)} \right\}\left(\frac{\xi(x)}{4x-1}\right)^\frac14 \Ai(-\xi(x))
$$
as $q\to 1^-$.

\section{Numerical Verification}
In the following tables, we have used a Maple-aided program
to verify the asymptotic formulas in Theorem~\ref{thm:1} and the limiting relation in Theorem~\ref{thm:2}.
The values are represented in scientific notations; for example,
$$
-2.325\text{e-}3=-2.325\times10^{-3}=-0.00235.
$$
In Table 1, `True' stands for the true value of $(q;q)_nS_n(uq^{-nt};q)$ by summing the series in (\ref{eq:1.4});
`Approx.' stands for the approximate value obtained from (\ref{eq:3.13});
`Error' is the relative error of the approximate value.
The degree of the polynomial is $n=50$ and $q=0.5$ in Table 1,
$u$ takes the values of $1$, $-1$ and $1+i$,
and $t$ takes the values of $0,0.5,1.0,1.2$ and $1.6$.
Here, we would like to mention that the approximate values are very close to the true values.
For instance, in the case when $u=t=1$,
the true value is -5.83981318477869$\cdots\times10^{187}$,
whereas the approximate value is -5.83981318477868$\cdots\times10^{187}$.
If we take only a few digits as what we have done in Table 1,
the two values all appear to be nearly the same.
In Table 2, `True' stands for the true value of $A_q(\sqrt q x)$ by summing the series in (\ref{eq:8});
`Approx.' stands for the approximate value obtained from the quantity on the right-hand side of (\ref{eq:13});
`Error' is the relative error of the approximate value.
We examine the cases when $x=0.5,1.0,4.0,10,20$ and $q=0.9,0.92,0.94,0.96,0.98,0.99$.

\begin{table}[h]
\caption{Numerical Verification of Theorem 1}
\begin{tabular}{|c|c|cccccc|}
\hline
\backslashbox{$u$\kern-1em}{$t$}
         &        &0          &0.5          &0.8           &1.0            &1.2            &1.6          \\
\hline\multirow{3}*{1}
         &True    &0.16076    &-9.3534e42   &1.0831e120    &-5.8398e187   &3.5453e270      &1.8649e481  \\
         \cline{2-8}
         &Approx. &0.16076    &-9.3534e42   &1.0831e120    &-5.8398e187   &3.5453e270      &1.8649e481\\
         \cline{2-8}
         &Error   &2.99e-15   &2.98e-8      &1.78e-15      &1.18e-15      &6.05e-13        &6.36e-7\\
\hline\multirow{3}*{-1}
         &True    &2.17267    &8.0063e47    &1.9036e121    &1.0264e189    &6.2313e271      &3.2740e482   \\
         \cline{2-8}
         &Approx. &2.17267    &8.0063e47    &1.9306e121    &1.0264e189    &6.2312e271      &3.2778e482 \\
         \cline{2-8}
         &Error   &6.31e-16   &6.12e-12     &1.11e-9       &3.54e-8       &1.13e-6         &1.16e-3 \\
\hline
\end{tabular}
\begin{tabular}{|c|c|cccc|}
\hline
\backslashbox{$u$\kern-1em}{$t$}
         &        &0               &0.5               &  1.0                &1.6\\
\hline\multirow{3}*{1+i}
         &True    &0.0117-0.6786i  &(1.92+8.38i)e48   &(-8.18-1.87i)e191    &(4.107-2.571i)e487    \\
         \cline{2-6}
         &Approx. &0.0117-0.6786i  &(1.92+8.38i)e48   &(-8.18-1.87i)e191    &(4.106-2.578i)e487    \\
         \cline{2-6}
         &Error   &8.83e-17        &3.16e-12          &1.39e-8              &4.55e-4       \\
\hline
\end{tabular}
\end{table}

\begin{table}[h]
\caption{Numerical Verification of Theorem 2}
\begin{tabular}{|c|c|cccccc|}
\hline
\backslashbox{$x$\kern-1em}{$q$}
&              & 0.9             & 0.92         & 0.94        &0.96         &0.98          &0.99\\
\hline\multirow{3}*{0.5}
&True          &-2.325e-3        &-3.826e-4     &1.120e-5     &-2.966e-8    &5.080e-16     &4.9298e-32 \\
\cline{2-8}
&Approx.       &-2.320e-3        &-3.819e-4     &1.118e-5     &-2.964e-8    &5.078e-16     &4.9303e-32\\
\cline{2-8}
&Error         & 0.0022          &0.0018        & 0.0012      &0.00080      &0.00038       &0.0000995\\
\hline\multirow{3}*{1.0}
&True          &   -5.171e-4    &2.978e-5     &-2.556e-6     &1.326e-9    &2.178e-18  &4.1417e-36 \\
\cline{2-8}
&Approx.       &   -5.159e-4    &2.973e-5     &-2.553e-6     &1.325e-9    &2.177e-18  &4.1408e-36\\
\cline{2-8}
&Error         &     0.0022     &0.0018       & 0.0013       &0.00087    &0.00043    &0.00021\\
\hline\multirow{3}*{4.0}
&True           &0.034973         &0.01680       &4.4202e-4          &-1.0084e-4              &4.4912e-8        &-5.6869e-16 \\
\cline{2-8}
&Approx.        &0.034891         &0.01677        &4.1898e-4         &-1.0073e-4               &4.4893e-8        &-5.6853e-16 \\
\cline{2-8}
&Error          & 0.00235          &0.0018        & 0.0028             &0.00107                 &0.00043          & 0.00028\\
\hline\multirow{3}*{10}
&True                       &38.6522         &-247.876     &2715.83      &43744.8       &3.3978e10     &2.1941e21 \\
\cline{2-8}
&Approx.                  &38.5316          &-247.372     &2712.29      &43745.2       &3.3961e10     &2.1944e21\\
\cline{2-8}
&Error                 & 0.0031          &0.0020       & 0.0013      &0.00022     &0.00051        &0.00014\\
\hline\multirow{3}*{20}
&True                &-2.0951e5         &-3.5927e6     &-5.9716e9     &7.3472e14       &-6.1129e29     &1.5900e61 \\
\cline{2-8}
&Approx.           &-2.0884e5         &-3.5801e6     &-5.9645e9     &7.3418e14       &-6.1124e29     &1.5897e61\\
\cline{2-8}
&Error               & 0.00320          &0.00349       & 0.00119      &0.00073         &0.00007        &0.00023\\
\hline
\end{tabular}
\end{table}


\begin{thebibliography}{9}


\bibitem{BaikSuidan} J. Baik and T.M. Suidan,
Random matrix central limit theorems for nonintersecting random walks,
\textit{Ann. Probab.} \textbf{35} (2007), 1807-1834.

\bibitem{Chihara}T.S. Chihara,
\textit{An Introduction to Orthogonal Polynomials}, Gordon and Breach, New York, 1978.

\bibitem{Christiansen}J.S. Christiansen,
The moment problem associated with the Stieltjes-Wigert polynomials,
\textit{J. Math. Anal. Appl.} \textbf{277} (2003), 218-245.


\bibitem{DT}
Y. Dolivet and M. Tierz,
Chern-Simons matrix models and Stieltjes-Wigert polynomials,
\textit{J. Math. Phys.}, \textbf{48} (2007), Arcile ID: 0235207, 20pp.

\bibitem{Ismail} M.E.H. Ismail, \textit{Classical and Quantum Orthogonal Polynomials in One Variable}, Cambridge University Press, Cambridge, 2005.

\bibitem{Ismail2} M.E.H. Ismail, Asymptotics of $q$-orthogonal polynomials and a $q$-Airy function, \textit{Int. Math. Res. Not.} \textbf{18} (2005), 1063-1088.

\bibitem{IsmailZhang} M.E.H. Ismail and R.-M. Zhang, Chaotic and periodic asymptotics for $q$-orthogonal polynomials, \textit{Int. Math. Res. Not.} (2006), Article ID 83274, 33 pp.

\bibitem{KS} R. Koekoek and R.F. Swarttouw,
The Askey-scheme of hypergeometric orthogonal polynomials and its $q$-analogue, Report No. 98-17, TU-Delft, 1998.

\bibitem{NIST}F.W.J. Olver, D.W. Lozier, R.F. Boisvert and C.W. Clark,
\textit{NIST Handbook of Mathematical Functions}, Cambridge University Press, New York, 2010.

\bibitem{Ramanujan}S. Ramanujan, \textit{The Lost Notebook and Other Unpublished Papers, Introduction by G.E. Andrews}, Narosa, New Delhi, 1988.

\bibitem{Szego} G. Szeg\"o,
\textit{Orthogonal Polynomials}, fourth ed., Colloquium Publications, vol. 23, Amer. Math. Soc., Providence, RI, 1975.

\bibitem{WangXSWong1} X.S. Wang and R. Wong, Discrete analogues of Laplace's approximation, \textit{Asymptot. Anal.} \textbf{54} (2007), 165-180.

\bibitem{WangXSWong2} X.S. Wang and R. Wong, Uniform asymptotics of some $q$-orthogonal polynomials, \textit{J. Math. Anal. Appl.} \textbf{364} (2010), 79-87.

\bibitem{WangWong} Z. Wang and R. Wong, Uniform asymptotics of the Stielties-Wigert polynomials via the Riemann-Hilbert approach, \textit{J. Math. Pures Appl.} \textbf{85} (2006), 698-718.

\bibitem{Wigert}S. Wigert,
Sur les polyn\^oes orthogonaux et l'approximation des fonctions continues,
\textit{Arkiv f\"or matematik, astronomi och fysik}, \textbf{17} (1923), 15pp.

\bibitem{Zhang} R. Zhang,
On a limiting relation between Ramanujan's entire function $A_q(z)$ and $\theta$-functions,
\textit{Q. J. Math.} \textbf{58} (2007), 519-532.

\bigskip
\bigskip

\end{thebibliography}
\end{document}